\documentclass[reqno]{amsart}

\usepackage{amsmath,amsthm,amsfonts,amscd,amssymb}
\usepackage[pdfborder={0 0 0}]{hyperref}
\usepackage{enumerate,verbatim}
\usepackage{enumitem}
\usepackage{cite}

\newcommand{\ZZ}{\mathbb{Z}}

\newcommand{\rank}{\operatorname{rank}}

\newtheorem{thm}{Theorem}[section]
\newtheorem{cor}[thm]{Corollary}
\newtheorem{lem}[thm]{Lemma}

\newtheorem{quest}{Question}

\theoremstyle{definition}

\theoremstyle{remark}
\newtheorem{rem}{Remark}[section]

\begin{document}

\title{Permutation invariant lattices}
\author{Lenny Fukshansky}
\author{Stephan Ramon Garcia}
\author{Xun Sun}\thanks{Fukshansky acknowledges support by Simons Foundation grant \#279155 and by the NSA grant \#1210223, Garcia acknowledges support by NSF grant DMS-1265973.}

\address{Department of Mathematics, 850 Columbia Avenue, Claremont McKenna College, Claremont, CA 91711, USA}
\email{lenny@cmc.edu}
\address{Department of Mathematics, Pomona College, 610 N. College Ave, Claremont, CA 91711, USA}
\email{stephan.garcia@pomona.edu}
\urladdr{\url{http://pages.pomona.edu/~sg064747/}}
\address{School of Mathematical Sciences, Claremont Graduate University, Claremont, CA 91711, USA}
\email{foxfur\_32@hotmail.com}


\begin{abstract} We say that a Euclidean lattice in~$\mathbb R^n$ is {\it permutation invariant} if its automorphism group has non-trivial intersection with the symmetric group~$S_n$, i.e., if the lattice is closed under the action of some non-identity elements of~$S_n$. Given a fixed element $\tau \in S_n$, we study properties of the set of all lattices closed under the action of~$\tau$: we call such lattices {\it $\tau$-invariant}. These lattices naturally generalize cyclic lattices introduced by Micciancio in~\cite{mic1, mic}, which we previously studied in~\cite{me_xun}. Continuing our investigation, we discuss some basic properties of permutation invariant lattices, in particular proving that the subset of well-rounded lattices in the set of all $\tau$-invariant lattices in~$\mathbb R^n$ has positive co-dimension (and hence comprises zero proportion) for all~$\tau$ different from an $n$-cycle.
\end{abstract}

\maketitle

\def\A{{\mathcal A}}
\def\AA{{\mathfrak A}}
\def\B{{\mathcal B}}
\def\C{{\mathcal C}}
\def\D{{\mathcal D}}
\def\EE{{\mathfrak E}}
\def\F{{\mathcal F}}
\def\G{{\mathcal G}}
\def\x{{\mathcal H}}
\def\I{{\mathcal I}}
\def\II{{\mathfrak I}}
\def\J{{\mathcal J}}
\def\K{{\mathcal K}}
\def\kk{{\mathfrak K}}
\def\L{{\mathcal L}}
\def\LL{{\mathfrak L}}
\def\M{{\mathcal M}}
\def\mm{{\mathfrak m}}
\def\MM{{\mathfrak M}}
\def\N{{\mathcal N}}
\def\O{{\mathcal O}}
\def\OO{{\mathfrak O}}
\def\PP{{\mathfrak P}}
\def\R{{\mathcal R}}
\def\W{{\mathcal W}}
\def\PNR{{\mathcal P_N(\real)}}
\def\PMNR{{\mathcal P^M_N(\real)}}
\def\PdNR{{\mathcal P^d_N(\real)}}
\def\s{{\mathcal S}}
\def\V{{\mathcal V}}
\def\X{{\mathcal X}}
\def\Y{{\mathcal Y}}
\def\Z{{\mathcal Z}}
\def\H{{\mathcal H}}
\def\cee{{\mathbb C}}
\def\Nn{{\mathbb N}}
\def\pee{{\mathbb P}}
\def\que{{\mathbb Q}}
\def\QQ{{\mathbb Q}}
\def\real{{\mathbb R}}
\def\RR{{\mathbb R}}
\def\zed{{\mathbb Z}}
\def\ZZ{{\mathbb Z}}
\def\aaa{{\mathbb A}}
\def\ff{{\mathbb F}}
\def\HDelta{{\it \Delta}}
\def\kk{{\mathfrak K}}
\def\qbar{{\overline{\mathbb Q}}}
\def\kbar{{\overline{K}}}
\def\ybar{{\overline{Y}}}
\def\kkbar{{\overline{\mathfrak K}}}
\def\ubar{{\overline{U}}}
\def\eps{{\varepsilon}}
\def\ahat{{\hat \alpha}}
\def\bhat{{\hat \beta}}
\def\k{{\nu}}
\def\gt{{\tilde \gamma}}
\def\h{{\tfrac12}}
\def\be{{\boldsymbol e}}
\def\bei{{\boldsymbol e_i}}
\def\bc{{\boldsymbol c}}
\def\bm{{\boldsymbol m}}
\def\bk{{\boldsymbol k}}
\def\bi{{\boldsymbol i}}
\def\bl{{\boldsymbol l}}
\def\bq{{\boldsymbol q}}
\def\bu{{\boldsymbol u}}
\def\bt{{\boldsymbol t}}
\def\bs{{\boldsymbol s}}
\def\bv{{\boldsymbol v}}
\def\bw{{\boldsymbol w}}
\def\bx{{\boldsymbol x}}
\def\bbx{{\overline{\boldsymbol x}}}
\def\bX{{\boldsymbol X}}
\def\bz{{\boldsymbol z}}
\def\bwy{{\boldsymbol y}}
\def\bY{{\boldsymbol Y}}
\def\bL{{\boldsymbol L}}
\def\ba{{\boldsymbol a}}
\def\bb{{\boldsymbol b}}
\def\bet{{\boldsymbol\eta}}
\def\bxi{{\boldsymbol\xi}}
\def\bo{{\boldsymbol 0}}
\def\bone{{\boldsymbol 1}}
\def\bol{{\boldsymbol 1}_L}
\def\ep{\varepsilon}
\def\p{\boldsymbol\varphi}
\def\q{\boldsymbol\psi}
\def\rank{\operatorname{rank}}
\def\aut{\operatorname{Aut}}
\def\lcm{\operatorname{lcm}}
\def\sgn{\operatorname{sgn}}
\def\spn{\operatorname{span}}
\def\md{\operatorname{mod}}
\def\Norm{\operatorname{Norm}}
\def\dim{\operatorname{dim}}
\def\det{\operatorname{det}}
\def\Vol{\operatorname{Vol}}
\def\rk{\operatorname{rk}}
\def\ord{\operatorname{ord}}
\def\ker{\operatorname{ker}}
\def\div{\operatorname{div}}
\def\Gal{\operatorname{Gal}}
\def\GL{\operatorname{GL}}
\def\SNR{\operatorname{SNR}}
\def\WR{\operatorname{WR}}
\def\IWR{\operatorname{IWR}}
\def\scg{\operatorname{\left< \Gamma \right>}}
\def\swrh{\operatorname{Sim_{WR}(\Lambda_h)}}
\def\ch{\operatorname{C_h}}
\def\cht{\operatorname{C_h(\theta)}}
\def\scgt{\operatorname{\left< \Gamma_{\theta} \right>}}
\def\scgmn{\operatorname{\left< \Gamma_{m,n} \right>}}
\def\gat{\operatorname{\Omega_{\theta}}}
\def\Obar{\operatorname{\overline{\Omega}}}
\def\Lbar{\operatorname{\overline{\Lambda}}}
\def\mn{\operatorname{mn}}
\def\disc{\operatorname{disc}}
\def\rot{\operatorname{rot}}
\def\Prob{\operatorname{Prob}}
\def\co{\operatorname{co}}
\def\ot{\operatorname{o_{\tau}}}
\def\oci{\operatorname{o}}
\def\Aut{\operatorname{Aut}}
\def\Mat{\operatorname{Mat}}
\def\SL{\operatorname{SL}}

\section{Introduction}
\label{intro}

Given a lattice $\Lambda \subset \real^n$ of rank $n$, $n \geq 2$, its successive minima
$$0 < \lambda_1 \leq \dots \leq \lambda_n$$
are defined as
$$\lambda_i = \min \left\{ r \in \real_{>0} : \dim_{\real} \spn_{\real} \left( \Lambda \cap B_n(r) \right) \geq i \right\},$$
where $B_n(r)$ is a ball of radius~$r$ centered at $\bo$ in $\real^n$. A collection of linearly independent vectors $\bx_1,\dots,\bx_n \in \Lambda$ such that $\|\bx_i\| = \lambda_i$, where $\|\cdot \|$ is the usual Euclidean norm, is referred to as a collection of vectors {\it corresponding to successive minima}; this collection is not unique, but there are only finitely many of them in a lattice with an upper bound on their number depending only on $n$. These vectors do not necessarily form a basis for $\Lambda$, however their span over $\zed$ is a sublattice of finite index in $\Lambda$. If $\Lambda$ has a basis consisting of vectors corresponding to successive minima, it is called a {\it Minkowskian lattice}. Vectors with norm equal to $\lambda_1$ are called {\it minimal vectors}. A lattice is called {\it well-rounded} (abbreviated WR) if $\lambda_1 = \dots = \lambda_n$, which is equivalent to saying that it has $n$ linearly independent minimal vectors. WR lattices are vital in extremal lattice theory and discrete optimization, as well as in connection with several other areas of mathematics; see~\cite{martinet} for a detailed overview. This provides motivation for studying distribution properties of WR lattices.

Let $S_n$ be the group of permutations on $n \geq 2$ elements. We can define an action of $S_n$ on $\real^n$ by $\tau(\bx) := (x_{\tau(1)},\dots,x_{\tau(n)})^t \in \real^n$ for each $\tau \in S_n$ and $\bx = (x_1,\dots,x_n)^t \in \real^n$. We say that a lattice $\Lambda \subset \real^n$ is {\it $\tau$-invariant} (or invariant under $\tau$) for a fixed $\tau \in S_n$ if $\tau(\Lambda) = \Lambda$. An important example of such lattices are lattices of the form 
\begin{equation}
\label{Ltx}
\Lambda_{\tau}(\bx) := \spn_{\zed} \left\{ \bx,\tau(\bx), \dots, \tau^{\nu-1}(\bx) \right\},
\end{equation}
where $\bx \in \real^n$ and $\nu$ is the order of $\tau$ in $S_n$. Cyclic lattices introduced by Micciancio in~\cite{mic1, mic} are precisely the full-rank sublattices of $\zed^n$ invariant under the $n$-cycle $\sigma_n := (1\ 2\dots n)$. In~\cite{me_xun} we investigated geometric properties of cyclic lattices, in particular studying well-rounded lattices of the form $\Lambda_{\sigma_n}(\bx)$ in every dimension. The goal of this note is to extend this investigation to more general permutation invariant lattices, as well as to outline some basic properties of permutation invariance.

There is another way to view the notion of permutation invariance. Given a lattice $\Lambda \subset \real^n$ of full rank, its automorphism group is defined as
$$\Aut(\Lambda) = \left\{ \sigma \in \GL_n(\zed) : \sigma(\bx) \cdot \sigma(\bwy) = \bx \cdot \bwy,\ \forall\, \bx,\bwy \in \Lambda \right\}.$$
It is a well known fact that $\Aut(\Lambda)$ is always a finite group; on the other hand, any finite subgroup of $\GL_n(\zed)$ is contained in the automorphism group of some lattice. In all dimensions except for $n=2,4,6,7,8,9,10$ (dimensions with exceptionally symmetric lattices) the largest such group is $(\zed/2\zed)^n \rtimes S_n$, the automorphism group of the integer lattice~$\zed^n$ (see~\cite{martinet} for more information). Lattices with large automorphism groups usually have a large degree of geometric symmetry. This often correlates with having many minimal vectors, and hence they have increased chances of being well-rounded. We show, however, that chances of being well-rounded are determined not just by the size of the automorphism group but by the type of elements it contains. Let $\tau \in S_n$ and define $\C_n(\tau)$ to be the set of all full rank $\tau$-invariant lattices in~$\real^n$. Then $\Lambda \in \C_n(\tau)$ if and only if $\tau \in \Aut(\Lambda)$.

\begin{quest} \label{WR-tau} Given $\tau \in S_n$, how big is the subset of WR lattices among all lattices in $\C_n(\tau)$?
\end{quest}

In order to make this question precise, we need some further notation. The space of all lattices in~$\real^n$ can be identified with the quotient space $\GL_n(\real)/\GL_n(\zed)$ via their basis matrices. Hence lattices are identified with points in a full dimensional subset of the $n^2$-dimensional Euclidean space $\Mat_{nn}(\real) \setminus \{A \in \Mat_{nn}(\real) : \det(A) = 0 \}$ modulo $\GL_n(\zed)$-equivalence. Each lattice has a Minkowski-reduced basis, consisting of short vectors, and reduction conditions amount to a finite collection of polynomial inequalities cutting out a full dimensional subset in the above set, called the {\it Minkowski reduction domain} (see~\cite{achill_book} for detailed information on its construction). WR lattices correspond to points in this space with $\lambda_1 = \dots = \lambda_n$, which can be interpreted as a collection of nontrivial polynomial conditions on the Minkowski reduction domain, and hence defines a subset of positive co-dimension. One can similarly talk about the dimension of the set $\C_n(\tau)$ by viewing it as a space of matrices with appropriate restrictions. With this notation in mind, we can state our main result.

\begin{thm} \label{main} Let $\tau \in S_n$ be different from an $n$-cycle. Then the subset of WR lattices in the set $\C_n(\tau)$ of full-rank $\tau$-invariant lattices in $\real^n$ has positive co-dimension.
\end{thm}

\noindent
In other words, the main difference between $n$-cycles and other permutations in this context is the following. If $\tau$ is an $n$-cycle, then lattices of the form $\Lambda_{\tau}(\bx)$ are {\it almost always} of full rank and are often well-rounded. On the other hand, when $\tau$ is not an $n$-cycle, lattices of the form $\Lambda_{\tau}(\bx)$ are {\it never} of full rank. 

This note is structured as follows. In Section~\ref{order} we prove that no vector can have $n$ linearly independent permutations by an element $\tau \in S_n$ unless $\tau$ is an $n$-cycle, discuss invariant subspaces of a permutation, and obtain an explicit formula for the number of linearly independent permutations of a generic vector.  We then use these observations in our proof of Theorem~\ref{main}, which is presented in Section~\ref{WR_prop}, along with a more geometric interpretation of this result.

The proof of Theorem~\ref{main} is based on showing that the polynomial conditions defining the subset of WR lattices remain nontrivial on $\C_n(\tau)$ whenever $\tau$ is not an $n$-cycle, and hence carve out a subset of $\C_n(\tau)$ of positive co-dimension. On the other hand, these conditions are automatically satisfied (i.e., become trivial) for a positive proportion of lattices of the form $\Lambda_{\tau}(\bx)$ when $\tau$ is an $n$-cycle, as discussed in~\cite{me_xun}. We are now ready to proceed.
\bigskip

\section{Permutation invariance properties}
\label{order}

In this section we discuss several basic properties of permutation invariance with respect to a fixed element of~$S_n$. Let $\tau \in S_n$ be an element of order $\k$, and for each vector $\bx \in \real^n$, define the $n \times \k$ matrix
$$M_{\tau}(\bx) := \left( \bx\ \tau(\bx)\ \dots\ \tau^{\k-1}(\bx) \right).$$
Then 
$$\Lambda_{\tau}(\bx) := M_{\tau}(\bx) \zed^\k = \spn_{\zed} \left\{ \bx,\tau(\bx),\dots,\tau^{\k-1}(\bx) \right\}$$
is  a lattice, as in~\eqref{Ltx} above, and we define the {\it $\tau$-order} of $\bx$ to be $\ot(\bx) := \rk \left( \Lambda_{\tau}(\bx) \right)$. We start by giving a bound on~$\ot(\bx)$ for an arbitrary vector in~$\real^n$. Let $\tau$ have a decomposition into disjoint cycles
\begin{equation}
\label{decomp}
\tau = c_1 \dots c_{\ell},
\end{equation}
where each $c_i \in S_n$ is a $k_i$-cycle and $\sum_{i=1}^{\ell} k_i = n$. The following lemma will serve as an important tool in the proof of our main result in Section~\ref{WR_prop}.

\begin{lem} \label{t_order} Let $\bx \in \real^n$, then
\begin{equation}
\label{ot-form}
\ot(\bx) \leq n - (\ell-1).
\end{equation}
\end{lem}

\proof
Let $T$ be the permutation matrix representation of $\tau$. Since $T$ is diagonalizable, the degree of its minimal polynomial $m_T(t)$ equals the number of distinct eigenvalues of $T$ (see, for instance, Corollary 3.3.10 of \cite{horn_johnson}). Let $\bx \in \real^n$, then
\begin{eqnarray}
\label{ot-bnd}
\ot(\bx) & = & \rk \Lambda_{\tau}(\bx) = \dim \spn_{\real} \left\{ \bx, \tau(\bx), \dots, \tau^{\nu-1}(\bx) \right\} \nonumber \\
& \leq & \dim \spn_{\real} \left\{ I,T,T^2,\dots,T^{\nu} \right\} \leq \deg m_T.
\end{eqnarray}
Let $\chi_T(t)$ be the characteristic polynomial of $T$. It is a well-known fact (see, for instance, Section~2.1 of~\cite{perm_matrix}) that
\begin{equation}
\label{char}
\chi_T(t) = \prod_{i=1}^{\ell} (t^{k_i}-1) =  \prod_{i=1}^{\ell} \left( \prod_{d \mid k_i} \Phi_d(t) \right),
\end{equation}
where $\Phi_d(t)$ is the irreducible $d$-th cyclotomic polynomial, whose roots are the $d$-th primitive roots of unity and whose degree is~$\varphi(d)$. The minimal and characteristic polynomials of $T$ have the same roots, although possibly with different multiplicities. Since  minimal polynomial of a permutation matrix $T$ must have distinct roots, these are precisely the distinct roots of~$\chi_T(t)$. Therefore to compute $\deg m_T$ we need to determine the number of distinct roots of~$\chi_T(t)$. Now~\eqref{char} implies that this number is equal to the sum of degrees of distinct cyclotomic polynomials dividing the product 
$$\prod_{i=1}^{\ell} (t^{k_i}-1),$$
and hence~\eqref{ot-bnd} combined with~\eqref{char} implies that
\begin{equation}
\label{ot-bnd-1}
\ot(\bx) \leq n - \sum_{ \substack{d \mid k_i, d \mid k_j \\ i < j}} \varphi(d) +\sum_{i<j<h}^{\ell} \sum_{d \mid \gcd(k_i,k_j,k_h)}\varphi(d) - \dots
\end{equation}
for any $\bx \in \real^n$, where dots in the upper bound of~\eqref{ot-bnd-1} indicate continued application of the inclusion-exclusion principle to obtain an actual formula for~$\deg m_T$; the sums above are understood as 0 if $\ell = 1$. 

On a simpler note, observe that the permutation matrix corresponding to a cycle of length $k$ contributes $k$ eigenvalues, one of which is necessarily $1$. Hence each cycle $c_i$ in the product~\eqref{decomp} for $i > 1$ actually contributes at most $k_i-1$ distinct eigenvalues, so
$$\ot(\bx) \leq n-(\ell-1)$$
for any $\bx \in \real^n$. This completes the proof of the lemma.
\endproof

\begin{rem} \label{ot<n} In particular, \eqref{ot-bnd-1} implies that unless $\tau$ is an $n$-cycle, $\ot(\bx) < n$ for any~$\bx \in \real^n$.
\end{rem}
\smallskip

Next we discuss invariant subspaces of a given permutation. We say that a proper subspace $V \subset \real^n$ is $\tau$-invariant if $\tau(V)=V$. We make two brief remarks on invariant subspaces of permutations. 

\begin{lem} \label{inv_sub} Suppose $\tau \in S_n$ is not an $n$-cycle. Then there exist infinitely many $\tau$-invariant subspaces of $\real^n$.
\end{lem}

\proof
Suppose $\tau$ is a product of $\ell$ disjoint cycles, $\ell \geq 2$. Pick any real numbers $a_1,\dots,a_\ell$, and define a vector $\bx \in \real^n$ by setting $x_i=a_j$ if $i$ is in the $j$-th cycle. Then $\spn_{\real}\{ \bx \}$ is a $\tau$-invariant subspace.  Since $\tau$ is not an $n$-cycle, there clearly are
infinitely many such subspaces.
\endproof

The situation is different when $\tau$ is an $n$-cycle, as we discussed in~\cite{me_xun}. Let $\sigma_n$ be the $n$-cycle $(1 \dots n)$, then for any $n$-cycle $\tau \in S_n$ there exists $\rho \in S_n$ such that $\tau = \rho \sigma_n \rho^{-1}$. Given a vector $\ba \in \real^n$, write $\ba(t) = \sum_{i=1}^n a_i t^{i-1} \in \real[t]$. For each non-constant polynomial divisor $f(t)$ of $t^n-1$, define
\begin{equation}
\label{c1}
V(f) = \left\{ \ba \in \real^n : f(t) \mid \ba(t) \right\},
\end{equation}
which is a proper subspace of $\real^n$.

\begin{lem} \label{cyclic} Let $\tau = \rho \sigma_n \rho^{-1} \in S_n$ be an $n$-cycle. Then $W \subset \real^n$ is a $\tau$-invariant subspace of $\real^n$ if and only if $W=\rho(V(f))$ for some $V(f)$ as in~\eqref{c1}. This is a finite collection of proper subspaces of~$\real^n$, and $\ot(\ba) < n$ if and only if $\ba$ is in the union of these subspaces.
\end{lem}

\proof
It is an easy observation that $W$ is a $\tau$-invariant subspace if and only if $\rho^{-1}(W)$ is a $\sigma_n$-invariant subspace. Hence it is sufficient to prove that $\sigma_n$-invariant subspaces of $\real^n$ are precisely those of the form $V(f)$ for some $f(t) \mid t^n-1$, as in~\eqref{c1}. The fact that all such subspaces are~$\sigma_n$-invariant follows from Lemmas~2.3 and~2.4 of~\cite{peikert}. Now suppose that a proper subspace~$U$ of~$\real^n$ is $\sigma_n$-invariant.  Suppose that $U$ is not of the form $V(f)$ for any non-constant $f(t) \mid t^n-1$. Then there exists $\ba \in U$ such that $\ba(t)$ is not divisible by any non-constant factor of $t^n-1$, which implies that the collection of vectors $\ba,\sigma_n(\ba),\dots,\sigma_n^{n-1}(\ba)$ is linearly independent (see, for instance,~\cite{horn_johnson}), contradicting the assumption that $U$ is a proper subspace of~$\real^n$. 

The argument above in particular shows that $\ot(\ba) < n$ for any $\ba$ in the union of the finitely many $\tau$-invariant subspaces of~$\real^n$. On the other hand, suppose that $\ot(\ba) < n$, and let $W = \spn_{\real} \left\{ \ba,\tau(\ba),\dots,\tau^{n-1}(\ba) \right\}$. Then $W$ is a proper $\tau$-invariant subspace of~$\real^n$, and hence must be one of those described above.
\endproof

We say that a vector $\bx \in \real^n$ is generic (with respect to $\tau$) if $\ot(\bx)$ is as large as possible. We can now show that inequality~\eqref{ot-bnd-1} in the proof of Lemma~\ref{t_order} for the $\tau$-order becomes equality on generic vectors in~$\real^n$.

\begin{cor} \label{formula} Let $\tau \in S_n$ have a disjoint cycle decomposition as in~\eqref{decomp} above. For a generic vector $\bx \in \real^n$, there is equality in~\eqref{ot-bnd-1}.
\end{cor}

\proof
Let $\bx \in \real^n$, and let us write
$$\bx = (\bx_1\ \dots\ \bx_{\ell})^t,$$
where $\bx_i \in \real^{k_i}$ for each $1 \leq i \leq \ell$. Then
$$\tau(\bx) = (c_1(\bx_1)\ \dots\ c_\ell(\bx_\ell))^t,$$
and $\bx$ is generic if and only if none of the $\bx_i$'s are contained in any of the invariant subspaces of the cycles $c_1,\dots,c_\ell$. Then Lemma~\ref{cyclic} implies $\oci_{c_i}(\bx_i) = k_i$ for each $1 \leq i \leq \ell$, and so $\ot(\bx)$ must be equal to $n = \sum_{i=1}^\ell k_i$ minus the number of repeated eigenvalues of the permutation matrix of~$\tau$, i.e., the right hand side of~\eqref{ot-bnd-1}.
\endproof
\bigskip

\section{Proportion of well-rounded lattices}
\label{WR_prop}

In this section we prove our main result on the proportion of WR lattices among permutation invariant lattices. As above, let $\tau \in S_n$ be an element of order $\k$ and let $\C_n(\tau)$ be the set of all full-rank $\tau$-invariant lattices in~$\real^n$. 

Let $\G_n(\tau)$ be the subset of all Minkowskian lattices in $\C_n(\tau)$. Naturally, every lattice $\Gamma \in \C_n(\tau)$ has a Minkowskian sublattice $\Gamma' \in \G_n(\tau)$ spanned by the vectors corresponding to successive minima of $\Gamma$. While Minkowskian sublattices may not be unique, there can only be finitely many of them, where an upper bound on this number depends only on~$n$. On the other hand, the index $|\Gamma : \Gamma'|$ of a Minkowskian sublattice is also bounded from above by a constant depending only on~$n$, and hence a given lattice in $\G_n(\tau)$ can be a Minkowskian sublattice for only finitely many lattices in~$\C_n(\tau)$ (see~\cite{martinet1} and subsequent works of J. Martinet and his co-authors for more information on the index of Minkowskian sublattices). Furthermore, it is clear that a lattice $\Gamma \in \C_n(\tau)$ is WR if and only if each corresponding lattice $\Gamma' \in \G_n(\tau)$ is WR.

As discussed in Section~\ref{intro} above, our goal in this note is to address the question of how often lattices in $\C_n(\tau)$ are well-rounded. From the discussion above, it is evident that we can instead consider lattices in $\G_n(\tau)$ and ask how often are they WR. We first make some dimensional observations. If $\tau_1$ and $\tau_2$ are conjugate in $S_n$, then there exists $\rho \in S_n$ such that $\tau_2=\rho \tau_1 \rho^{-1}$. Then $\Gamma \in \G_n(\tau_1)$ if and only if $\rho(\Gamma) \in \G_n(\tau_2)$ and $\Gamma$ is WR if and only $\rho(\Gamma)$ is WR. Furthermore, the sets $\G_n(\tau_1)$ and $\G_n(\tau_2)$ have the same dimension. Define
$$\W_n(\tau) = \{ \Gamma \in \G_n(\tau) : \Gamma \text{ is WR} \}.$$

\begin{lem} \label{cyclic_dim} Let $\tau \in S_n$ be an $n$-cycle, define
$$\C'_n(\tau) = \left\{ \Gamma \in \C_n(\tau) : \Gamma = \Lambda_{\tau}(\bx) \text{ for some } \bx \in \real^n \right\},$$
$$\G'_n(\tau) = \left\{ \Gamma \in \G_n(\tau) : \Gamma = \Lambda_{\tau}(\bx) \text{ for some } \bx \in \real^n \right\},$$
$$\W'_n(\tau) = \left\{ \Gamma \in \W_n(\tau) : \Gamma = \Lambda_{\tau}(\bx) \text{ for some } \bx \in \real^n \right\}.$$
Then $n \leq \dim \W'_n(\tau) \leq \dim \G'_n(\tau) = \dim \C'_n(\tau)$.
\end{lem}

\proof
By the remark above, we can assume without loss of generality that $\tau = \sigma_n := (1 \dots n)$, and it is also clear that $\C_n(\tau)$ and $\G_n(\tau)$ have the same dimension. Then the statement of the lemma follows from Lemma~3.2 of~\cite{me_xun}.
\endproof

\begin{lem} \label{dim} Let $\tau = c_1 \dots c_{\ell}$ be a disjoint cycle decomposition for $\tau$ so that $c_i$ is an $m_i$-cycle ($m_i > 1$) for every $1 \leq i \leq \ell$ and $q=n-\sum_{i=1}^{\ell} m_i$. Then the dimension of the sets $\C_n(\tau)$ and $\G_n(\tau)$ is at least $n+q(q-1)$. This includes the case of $\tau$ being the identity in $S_n$, when we have $\ell=0$ and $n=q$.
\end{lem}

\proof
For each $1 \leq i \leq \ell$, let $\Lambda_i \subset \real^{m_i}$ be a lattice in $\C_{m_i}(\sigma_{m_i})$, where $\sigma_{m_i}$ is the standard $m_i$-cycle $(1 \dots m_i)$, and let $\Omega \subset \real^q$ be any lattice in $\real^q$. Let us view each of the lattices $\Lambda_i$ and $\Omega$ as embedded into $\real^n$: each $\Lambda_i$ into the $\sum_{j=1}^{i-1} m_j + 1$ through $\sum_{j=1}^i m_i$ coordinates of $\real^n$, and $\Omega$ into the last $q$ coordinates. Then the lattice
\begin{equation}
\label{dim1}
\Lambda = \Lambda_1 \oplus \dots \oplus \Lambda_{\ell} \oplus \Omega \subset \real^n
\end{equation}
is a lattice in $\C_n(\tau)$. By Lemma~\ref{cyclic_dim}, each $\C_{m_i}(\sigma_{m_i})$ has dimension at least $m_i$, while the dimension of the space of all lattices in $\real^q$ is $q^2$. Hence dimension of the set of all lattices as in~\eqref{dim1} is
$$\sum_{i=1}^{\ell} m_i + q^2 = n - q + q^2,$$
which completes the proof of the lemma.
\endproof

Next we consider the dimension of the set of WR lattices in~$\G_n(\tau)$ when $\tau$ is not an $n$-cycle.
 
\begin{lem} \label{WR_dim} Let $\tau \in S_n$ not be an $n$-cycle, then $\tau$ has a decomposition into $\ell > 1$ disjoint cycles, as in~\eqref{decomp}. Define $\ot$ to be the $\tau$-order of a generic vector in~$\real^n$ as above. Then
\begin{equation}
\label{dim_wr}
\dim \W_n(\tau) \leq \dim \G_n(\tau) - \left( \left \lceil{\frac{n}{\ot}}\right \rceil - 1 \right) < \dim \G_n(\tau).
\end{equation}
\end{lem}

\proof
A trivial implication of Lemma~\ref{dim} is that the set $\G_n(\tau)$ is not empty for any $\tau \in S_n$. Furthermore, the set $\W_n(\tau)$ is also not empty, since $\zed^n \in \W_n(\tau)$ for every $\tau \in S_n$. Let $\Lambda \in \G_n(\tau)$, then in order for $\Lambda$ to be in $\W_n(\tau)$ the following conditions have to be satisfied by the successive minima of $\Lambda$:
\begin{equation}
\label{sm_cond}
\lambda_1(\Lambda) = \dots = \lambda_n(\Lambda).
\end{equation}
We now show that there are at least $\left \lceil{\frac{n}{\ot}}\right \rceil - 1$ nontrivial polynomial conditions on $\G_n(\tau)$ among equations~\eqref{sm_cond}. Let us write $S(\Lambda)$ for the set of minimal vectors of $\Lambda$, then $\Lambda \in \W_n(\tau)$ if and only if $S(\Lambda)$ contains $n$ linearly independent vectors. Suppose that $\ba \in S(\Lambda)$, then $\tau^k(\ba) \in S(\Lambda)$ for each $k$ and automatically,
$$\|\ba\| = \|\tau^k(\ba)\|$$
for each $k$. From Lemma~\ref{t_order} we know that $\ot(\ba) \leq \ot < n$, hence there must exist some $\bb \in S(\Lambda)$ such that $\bb \neq \tau^k(\ba)$ for any $k$. On the other hand, we must still have
$$\|\ba\| = \|\bb\|,$$
which now constitutes a nontrivial equation. Hence the polynomial conditions~\eqref{sm_cond} defining $\W_n(\tau)$ as a subset of $\G_n(\tau)$ are not automatically satisfied, and thus carve out a subset of positive co-dimension. In fact, the number of such equations is at least the number of blocks of the form $\Lambda_{\tau}(\bx)$ required to span a lattice of rank~$n$ minus 1, which is~$\left \lceil{\frac{n}{\ot}}\right \rceil - 1 > 1$. This implies~\eqref{dim_wr}.
\endproof

\proof[Proof of Theorem~\ref{main}] The theorem now follows from Lemma~\ref{WR_dim} and our previous observation that for any lattice in~$\G_n(\tau)$ there are only finitely many in~$\C_n(\tau)$ and vice versa, where these proportions depend only on the dimension~$n$.
\endproof
\smallskip

There is also another geometric interpretation of our result. Notice that the set $\G_n(\tau)$ is closed under isometries, i.e., a lattice $\Gamma$ is isometric to some lattice $\Gamma' \in \G_n(\tau)$ only if $\Gamma$ itself is in $\G_n(\tau)$. It is a well known fact (see, for instance,~\cite{gruber} or~\cite{martinet}) that lattices up to isometry are in bijective correspondence with positive definite quadratic forms up to arithmetic equivalence. In other words, isometry classes of matrices in $\GL_n(\real)/\GL_n(\zed)$ are in bijective correspondence with arithmetic equivalence classes of positive definite $n \times n$ symmetric matrices modulo $\GL_n(\zed)$. Specifically, this correspondence is given by $A \to A^tA$ where $A$ is a basis matrix of a full-rank lattice in~$\real^n$. The conditions on the basis being Minkowski reduced translate into a finite collection of linear inequalities in the entries of $A^tA$, which define a polyhedral cone of full dimension with finitely many facets in the $\frac{n(n+1)}{2}$-dimensional space of $n \times n$ symmetric matrices; this cone is called the Minkowski reduction domain (see~\cite{achill_book} for details). Points in the interior of each facet of this cone correspond to lattices with the same automorphism group, i.e., automorphism groups are constant on the open facets and lattices in the interior of the cone have the trivial automorphism group~$\{\pm 1\}$; this fact goes back to the work of Minkowski. Furthermore, the set of all quadratic forms in $n$ variables whose automorphism groups contain a fixed subgroup of $\GL_n(\zed)$ is known to be a linear space, called a Bravais manifold (see, for instance,~\cite{achill_book, tammela}), and the positive definite forms (which correspond to lattices) in this space form a polyhedral cone. On the other hand, the subset of well-rounded lattices has co-dimension~$n$ in the space of all full-rank lattices in~$\real^n$ and is known to be a minimal $\SL_2(\zed)$-equivariant contractible deformation retract of this space (see Remark~3.3 of~\cite{lizhen_ji} for a detailed overview of this construction, which was originally introduced in the work of Soule, Lannes, and Ash; see also~\cite{pettet} for the proof of minimality of this retract). Then our Theorem~\ref{main} can be interpreted as a certain characterization of the ``size" of intersection of different facets of the Minkowski domain and Bravais manifolds with the well-rounded retract.
\bigskip

{\bf Acknowledgement.} We thank Bryan Ek for noticing a mistake in the original version of our Lemma~\ref{t_order} and suggesting a correction.
\bigskip

\bibliographystyle{plain}  
\bibliography{perm_invariant}    

\begin{thebibliography}{10}

\bibitem{me_xun}
L.~Fukshansky and X.~Sun.
\newblock On the geometry of cyclic lattices.
\newblock {\em Discrete Comput. Geom.}, 52(2):240--259, 2014.

\bibitem{gruber}
P.~Gruber and C.~G. Lekkerkerker.
\newblock {\em Geometry of Numbers}.
\newblock North-Holland Mathematical Library, Second Edition, 1987.

\bibitem{perm_matrix}
B.~M. Hambly, P.~Keevash, N.~O'Connell, and D.~Stark.
\newblock The characteristic polynomial of a random permutation matrix.
\newblock {\em Stochastic Process. Appl.}, 90(2):335--346, 2000.

\bibitem{horn_johnson}
R.~A. Horn and C.~R. Johnson.
\newblock {\em Matrix Analysis}.
\newblock Cambridge University Press, Cambridge, 2nd edition, 2013.

\bibitem{lizhen_ji}
L.~Ji.
\newblock Well-rounded equivariant deformation retracts of
  {T}eichm$\ddot{\mathrm{u}}$ller spaces.
\newblock {\em to appear in Enseign. Math.; arXiv:1302.0877}, 2014.

\bibitem{martinet1}
J.~Martinet.
\newblock Sur l'indice d'un sous-r\'eseau.
\newblock In {\em R\'eseaux euclidiens, designs sph\'eriques et formes
  modulaires, Monogr. Enseign. Math., 37}, pages 163--211. Enseignement Math.,
  Geneva, 2001.

\bibitem{martinet}
J.~Martinet.
\newblock {\em Perfect Lattices in Euclidean Spaces}.
\newblock Springer-Verlag, 2003.

\bibitem{mic1}
D.~Micciancio.
\newblock Generalized compact knapsacks, cyclic lattices, and efficient one-way
  functions from worst-case complexity assumptions.
\newblock {\em FOCS, IEEE Computer Society}, pages 356--365, 2002.

\bibitem{mic}
D.~Micciancio.
\newblock Generalized compact knapsacks, cyclic lattices, and efficient one-way
  functions.
\newblock {\em Comput. Complexity}, 16(4):365--411, 2007.

\bibitem{peikert}
C.~Peikert and A.~Rosen.
\newblock Efficient collision-resistant hashing from worst-case assumptions on
  cyclic lattices.
\newblock {\em Theory of cryptography, Lecture Notes in Comput. Sci., 3876,
  Springer, Berlin,}, pages 145--166, 2006.

\bibitem{pettet}
A.~Pettet and J.~Souto.
\newblock Minimality of the well-rounded retract.
\newblock {\em Geom. Topol.}, 12(3):1543--1556, 2008.

\bibitem{achill_book}
A.~Sch$\ddot{\mathrm{u}}$rmann.
\newblock {\em Computational geometry of positive definite quadratic forms},
  volume~48 of {\em University Lecture Series}.
\newblock American Mathematical Society, Providence, RI, 2009.

\bibitem{tammela}
P.~P. Tammela.
\newblock Two remarks on {B}ravais varieties.
\newblock {\em Investigations in number theory, 2. Zap. Nauchn. Sem. Leningrad.
  Otdel. Mat. Inst. Steklov. (LOMI)}, 3:90--93, 1973.

\end{thebibliography}
\end{document}